%% file: _main.tex
\documentclass[final,5p,times,twocolumn]{elsarticle}
\pdfoutput=1

\usepackage{amsmath}
\usepackage{amssymb}
\usepackage{amsthm}
\usepackage[short,nocomma]{optidef}
\usepackage{stackengine,scalerel}
\usepackage{multirow}
\usepackage{xfrac}
\usepackage{xcolor}
\usepackage{mathtools,bm}

\usepackage{titlesec}
\titleformat*{\section}{\bfseries\unboldmath}

\usepackage{hyperref}

\newcommand*\diff{\mathop{}\!\mathrm{d}}

\DeclareMathOperator{\prob}{Pr}

\def\bs{\boldsymbol}

\DeclareMathOperator{\tr}{tr}

\DeclarePairedDelimiterX\Set[1]\{\}{%
  #1%
}

\newtheorem{thm}{Theorem}
\newtheorem{cor}[thm]{Corollary}
\newtheorem{lem}[thm]{Lemma}

\theoremstyle{definition}
\newtheorem{defn}{Definition}
\newtheorem{example}{Example}
\theoremstyle{remark}

\newcommand{\Solset}{\mathcal{U}}   
\newcommand{\Sol}{U}  
\newcommand{\sol}{u}  
\DeclareRobustCommand{\Optsol}{\overline{U}}  
\newcommand{\optsol}{\overline{u}}  

\newcommand{\optorth}{\tilde{u}}
\newcommand{\Optorth}{\widetilde{U}}
\newcommand{\Orthset}{\widetilde{\mathcal{U}}}

\newcommand{\perm}{m}
\newcommand{\Perm}{M}
\newcommand{\Permset}{\mathcal{P}}

\newcommand{\Optinv}{\overline{V}}  
\newcommand{\optinv}{\overline{v}}  

\newcommand{\srcwt}{p}                  
\newcommand{\srcwtvec}{\boldsymbol{p}}  
\newcommand{\srctm}{S}                  
\newcommand{\srctmvec}{\boldsymbol{S}}  

\newcommand{\tgtwt}{q}                  
\newcommand{\tgtwtvec}{\boldsymbol{q}}  
\newcommand{\tgttm}{T}                  
\newcommand{\tgttmvec}{\boldsymbol{T}}  

\newcommand{\dualavec}{\boldsymbol{\alpha}}
\newcommand{\optduala}{\overline{\alpha}}
\newcommand{\optdualavec}{\overline{\boldsymbol{\alpha}}}

\newcommand{\dualbvec}{\boldsymbol{\beta}}
\newcommand{\optdualb}{\overline{\beta}}
\newcommand{\optdualbvec}{\overline{\boldsymbol{\beta}}}

\newcommand{\dualcmat}{\Gamma}
\newcommand{\optdualc}{\overline{\gamma}}
\newcommand{\optdualcmat}{\overline{\Gamma}}

\journal{Operations Research Letters}

\begin{document}

\begin{frontmatter}

\title{On Simplices with a Given Barycenter That Are Enclosed by the Standard Simplex}

\author{Brent Austgen\corref{cor1}}
\ead{brent.austgen@utexas.edu}
\author{John J. Hasenbein}
\author{Erhan Kutanoglu}
\cortext[cor1]{Corresponding author}
\address{Operations Research \& Industrial Engineering Program, The University of Texas at Austin, Austin, TX, USA}

\begin{abstract}
We present an optimization model defined on the manifold of the set of stochastic matrices. Geometrically, the model is akin to identifying a maximum-volume $n$-dimensional simplex that has a given barycenter and is enclosed by the $n$-dimensional standard simplex. Maximizing the volume of a simplex is equivalent to maximizing the determinant of its corresponding matrix. In our model, we employ trace maximization as a linear alternative to determinant maximization. We identify the analytical form of a solution to this model. We prove the solution is optimal and present necessary and sufficient conditions for it to be the unique optimal solution. Additionally, we show the identified optimal solution is an inverse $M$-matrix, and that its eigenvalues are the same as its diagonal entries. We demonstrate how the model and its solutions apply to the task of synthesizing conditional cumulative distribution functions (CDFs) that, in tandem with a given discrete marginal distribution, coherently preserve a given CDF.
\end{abstract}

\begin{highlights}
\item We present an optimization model that is akin to identifying a maximum-volume simplex that has a given barycenter and is enclosed by the standard simplex.
\item We provide the analytical form of a particular solution to the model, a matrix $\Optsol$, and prove that it is optimal. Additionally, we prove necessary and sufficient conditions for $\Optsol$ to be the unique optimal solution.
\item We provide the analytical form of $\Optinv = \Optsol^{-1}$. We prove that $\Optinv$ is an $M$-matrix and that $\Optsol$ is thus an inverse $M$-matrix. We also prove that the eigenvalues of $\Optsol$ and $\Optinv$ are the same as their diagonal entries.
\item We demonstrate how the optimization model and its solutions apply to the task of synthesizing conditional cumulative distribution functions (CDFs) that, in tandem with a given discrete marginal distribution, coherently preserve a given CDF.
\end{highlights}


\end{frontmatter}


\input{sections/01_introduction}
\input{sections/02_optimization}
\input{sections/03_solutions}
\input{sections/04_properties}
\input{sections/05_application}
\input{sections/06_conclusion}

\section*{Competing Interests}
The authors have no competing interests to declare.

\section*{Data Availability}
No data was used for the research described in the article.


\bibliographystyle{elsarticle-num-names}\biboptions{authoryear}
\bibliography{_main}

\end{document}

%% file: sections/01_introduction.tex
\section{Introduction}

\subsection{Overview}

We present an optimization model defined on the manifold of the set of stochastic matrices. Geometrically, the model is akin to identifying a maximum-volume $n$-dimensional simplex. The identification of maximum- and minimum-volume simplices under various conditions has many applications including but not limited to hyperspectral unmixing \citep{Craig1994,Li2008,Hendrix2013}, system control \citep{Lombardi2009,Kuti2014}, computer-aided design \citep{Tordesillas2022}, and signal processing \citep{Douik2019}. In our model, the conditions we consider are that the simplex has a given barycenter and that the simplex is enclosed by the $n$-dimensional standard simplex. Maximizing the volume of a simplex is equivalent to maximizing the determinant of its corresponding matrix. In our model, we employ the trace (the sum of eigenvalues) as a linear alternative to the determinant (the product of eigenvalues).

We introduce the model and its dual in Section~\ref{section:optimization}. In Section~\ref{section:solutions}, we identify the analytical form of primal and dual solutions and prove that the solutions are optimal using the Karush-Kuhn-Tucker (KKT) conditions. We additionally prove a necessary and sufficient condition for the identified optimal primal solution, a matrix, to be uniquely optimal. In Section~\ref{section:properties}, we present the analytical inverse of that and then classify the two matrices. We additionally identify and discuss the eigenvalues of these matrices. In Section~\ref{section:application}, we present the application that motivates the optimization model and its solutions. The application involves synthesizing conditional cumulative distribution functions (CDFs) that, in tandem with a given discrete marginal distribution, coherently preserve a given CDF. We present our conclusions in Section~\ref{section:conclusion}.

\subsection{Notation, Conventions, \& Definitions} \label{subsection:notation}

In this paper, all matrices and vectors are indexed using 1\nobreakdash-based indexing. All vectors are column vectors unless transposed, and all transposed vectors are row vectors. We let $\bs{0}$ and $\bs{1}$ be vectors of all zeros or ones, respectively, and we let $\bs{e}_i$ be the vector with a 1 in the $i\textsuperscript{th}$ position and zeros elsewhere. We define $[c]^+ = \max\{c, 0\}$ and $[c]^- = \min\{c, 0\}$. We additionally define $\mathbb{I}_{[\,\cdot\,]}$ as the indicator function that takes a value of one if the subscripted condition is true and zero otherwise.

Our analysis revolves around stochastic matrices and vectors, terms we now define.
\begin{defn}
    A square matrix $A \in \mathbb{R}^{n \times n}$ is a \textit{(right) stochastic matrix} if $A \bs{1} = \bs{1}$ and $A \ge 0$ elementwise.
\end{defn}
\begin{defn}
    A vector $\bs{a} \in \mathbb{R}^n$ is a \textit{stochastic vector} if $\bs{a}^\top \bs{1} = 1$ and $\bs{a} \ge \bs{0}$ elementwise.
\end{defn}

In Section~\ref{section:application}, we discuss the 1-Wasserstein distance, also known as the Earth-Mover's distance. We now define that distance as it appears in \citet{Rachev1998}.
\begin{defn}
The \emph{1-Wasserstein distance} between probability measure $\mu_1$ with associated CDF $F_1(x)$ and probability measure $\mu_2$ with associated CDF $F_2(x)$ is
\begin{equation*}
    W_1(\mu_1, \mu_2) = \int_{\mathbb{R}} | F_1(x) - F_2(x) | \diff x.
\end{equation*}
\end{defn}

%% file: sections/02_optimization.tex
\section{Optimization Model} \label{section:optimization}

Let $\Sol \in \mathbb{R}^{n \times n}$ be a stochastic matrix and $\tgtwtvec \in \mathbb{R}^n$ a stochastic vector. Note that $\Sol \ge 0$ and $\tgtwtvec \ge \bs{0}$ together imply that $\Sol^\top \tgtwtvec \ge \bs{0}$, and that $\bs{1}^\top \Sol^\top \tgtwtvec = \bs{1}^\top \tgtwtvec = 1$. That is, $\Sol^\top$ maps one stochastic vector to another.

Given two $n$-dimensional stochastic vectors ${\tgtwtvec > \bs{0}}$ and ${\srcwtvec > \bs{0}}$, we consider as a candidate any stochastic matrix $U$ whose transpose maps $\tgtwtvec$ to $\srcwtvec$. That is, our candidate set is
\begin{equation}
    \Solset = \{
        U \in \mathbb{R}^{n \times n} :
        U^\top \bs{q} = \bs{p},
        U \bs{1} = \bs{1},
        U \ge 0
    \}.
\end{equation}
Among these matrices, we seek a matrix whose trace is maximized. The problem we describe may be formulated as a linear program (LP):
\begin{equation}
    \max\{\tr(\Sol) : \Sol \in \Solset\}. \label{eq:P} \tag{P}
\end{equation}
Recall that the trace of a matrix is the sum of its diagonal entries, which is also equal to the sum of its eigenvalues. As such, \eqref{eq:P} is an eigenvalue optimization problem.

The geometry of \eqref{eq:P} may be interpreted as follows. Take the rows of $\Sol$ to be vertices of an $(n-1)$\nobreakdash-simplex and the entries of $\tgtwtvec$ to be weights at those vertices. The feasible space is the set of $(n-1)$\nobreakdash-simplices that have a $\tgtwtvec$-weighted barycenter of $\srcwtvec$ and are enclosed by the standard ${(n-1)}$\nobreakdash-simplex. The goal is to identify a simplex for which the sum of $L^1$ distances between the $i\textsuperscript{th}$ vertex and $\bs{e}_i$ is minimized.

To obtain the dual of \eqref{eq:P}, consider its Lagrangian:
\begin{equation*}
    \mathcal{L}(\Sol, \dualavec, \dualbvec, \dualcmat) =
    \mathrm{tr}(\Sol)
    - \dualavec^\top (\Sol^\top \tgtwtvec - \srcwtvec)
    - \dualbvec^\top (\Sol \bs{1} - \bs{1})
    + \mathrm{tr}(\dualcmat \Sol^\top).
\end{equation*}
From this, we see that
\begin{equation*}
\nabla_\Sol \mathcal{L}(\Sol, \dualavec, \dualbvec, \dualcmat)
= I
- \tgtwtvec \dualavec^\top
- \dualbvec \bs{1}^\top
+ \dualcmat.
\end{equation*}
Hence, the dual function is
\begin{equation}
g(\dualavec, \dualbvec, \dualcmat)
= \begin{cases}
\dualavec^\top \srcwtvec + \dualbvec^\top \bs{1}~&\text{if}~ \tgtwtvec \dualavec^\top + \dualbvec \bs{1} ^\top - \dualcmat = I, \dualcmat \ge 0, \\
\infty~&\text{otherwise}.
\label{eq:dual_function}
\end{cases}
\end{equation}

The dual LP, derived from \eqref{eq:dual_function}, is
\begin{mini}
    {}{\dualavec^\top \srcwtvec + \dualbvec^\top \bs{1}}{\label{eq:D} \tag{D}}{}
    \addConstraint{\tgtwtvec \dualavec^\top + \dualbvec \bs{1} ^\top - \dualcmat = I,}
    \addConstraint{\dualcmat \ge 0.}
\end{mini}

We summarize the KKT conditions \citep{Karush2014,Kuhn2014} for \eqref{eq:P} and \eqref{eq:D} in Table~\ref{tab:KKT}. Since \eqref{eq:P} and \eqref{eq:D} are LPs, the KKT conditions are necessary and sufficient for optimality.
\begin{table}
    \centering
    \caption{KKT Conditions for \eqref{eq:P} and \eqref{eq:D}.}
    \label{tab:KKT}
    \renewcommand*{\arraystretch}{1.1}
    \begin{tabular}{|l|l|}
    \hline
    Stationarity & $\tgtwtvec \dualavec^\top + \dualbvec \bs{1} ^\top - \dualcmat = I$ \\
    \hline
    \multirow{3}{*}{Primal Feasibility} & $\srcwtvec^\top = \tgtwtvec^\top \Sol$ \\
    & $\Sol \bs{1} = \bs{1}$ \\
    & $\Sol \ge 0$ \\
    \hline
    Dual Feasibility & $\dualcmat \ge 0$ \\
    \hline
    Complementary Slackness & $\tr(\Sol \dualcmat^\top) = 0$ \\
    \hline
    \end{tabular}
\end{table}

%% file: sections/03_solutions.tex
\section{Optimal Primal \& Dual Solutions} \label{section:solutions}

For stochastic vectors $\srcwtvec > \bs{0}$ and $\tgtwtvec > \bs{0}$, let $\Optsol$ be the corresponding matrix with entries
\begin{equation}
    \optsol_{ij} \triangleq \begin{cases}
        1 - [1 - \srcwt_i / \tgtwt_i]^+,
        & i = j, \\[0.25em]
        0,
        & i \ne j~\text{and}~\srcwtvec = \tgtwtvec, \\[0.25em]
        \displaystyle \frac{[1 - \srcwt_i / \tgtwt_i]^+ [\srcwt_j - \tgtwt_j]^+}
                           {\sum_{k=1}^{n} [\srcwt_k - \tgtwt_k]^+},
        & i \ne j~\text{and}~\srcwtvec \ne \tgtwtvec.
    \end{cases}
\end{equation}
Additionally, let $\optdualavec$ and $\optdualbvec$ be $n$-dimensional vectors and $\optdualcmat$ an $n \times n$ matrix with entries
\begin{align*}
    \optduala_j
    &\triangleq \mathbb{I}_{[\srcwt_j < \tgtwt_j]} / \tgtwt_j, \\
    \optdualb_i
    &\triangleq \mathbb{I}_{[\srcwt_i \ge \tgtwt_i]}, \\
    \optdualc_{ij}
    &\triangleq \mathbb{I}_{[\srcwt_i \ge \tgtwt_i]}
              + \mathbb{I}_{[\srcwt_j < \tgtwt_j]} \cdot \tgtwt_i / \tgtwt_j
              - \mathbb{I}_{[i = j]}.
\end{align*}
We present these values on a case-by-case basis in Table~\ref{tab:dual_solution}.
\begin{table*}
    \centering
    \renewcommand*{\arraystretch}{1.1}
    \caption{Values of $\optdualavec$, $\optdualbvec$, $\optdualcmat$, and related quantities on a case-by-case basis.}
    \label{tab:dual_solution}
    \begin{tabular}{|c|c|c||c|c|c|c|c|}
        \hline
        \multicolumn{3}{|c||}{Conditions}
        & $\optduala_j$
        & $\optdualb_i$
        & $\optdualc_{ij}$
        & $\left[\tgtwtvec \optdualavec^\top\right]_{ij}$
        & $\left[\optdualbvec \bs{1}^\top\right]_{ij}$
        \\ \hline

        $i  =  j$ & $\srcwt_i  <  \tgtwt_i$ & $\srcwt_j  <  \tgtwt_j$
        & $1 / \tgtwt_j$
        & $0$
        & $0$
        & $1$
        & $0$
        \\ \hline

        $i  =  j$ & $\srcwt_i \ge \tgtwt_i$ & $\srcwt_j \ge \tgtwt_j$
        & $0$
        & $1$
        & $0$
        & $0$
        & $1$
        \\ \hline
        
        $i \ne j$ & $\srcwt_i  <  \tgtwt_i$ & $\srcwt_j  <  \tgtwt_j$
        & $1 / \tgtwt_j$
        & $0$
        & $\tgtwt_j / \tgtwt_i$
        & $\tgtwt_j / \tgtwt_i$
        & $0$
        \\ \hline

        $i \ne j$ & $\srcwt_i  <  \tgtwt_i$ & $\srcwt_j \ge \tgtwt_j$
        & $0$
        & $0$
        & $0$
        & $0$
        & $0$
        \\ \hline

        $i \ne j$ & $\srcwt_i \ge \tgtwt_i$ & $\srcwt_j  <  \tgtwt_j$
        & $1 / \tgtwt_j$
        & $1$
        & $\tgtwt_j / \tgtwt_i + 1$
        & $\tgtwt_j / \tgtwt_i$
        & $1$
        \\ \hline

        $i \ne j$ & $\srcwt_i \ge \tgtwt_i$ & $\srcwt_j \ge \tgtwt_j$
        & $0$
        & $1$
        & $1$
        & $0$
        & $1$
        \\ \hline
    \end{tabular}
\end{table*}

We now work toward proving that $\Optsol$ and $(\optdualavec, \optdualbvec, \optdualcmat)$ are an optimal primal-dual pair.
\begin{lem} \label{lem:primal_feasible}
For any two $n$-dimensional stochastic vectors ${\srcwtvec > \bs{0}}$ and ${\tgtwtvec > \bs{0}}$, the corresponding matrix $\Optsol$ is feasible in \eqref{eq:P}.
\end{lem}
\begin{proof}
For $\Optsol$ to be feasible in \eqref{eq:P}, it must exist in $\Solset$. If $\srcwtvec = \tgtwtvec$, then $\Optsol = I \in \Solset$, trivially. Otherwise, consider as preliminaries that
\begin{equation*}
    [1 - \srcwt_i / \tgtwt_i]^+ = [\tgtwt_i - \srcwt_i]^+ / \tgtwt_i,
\end{equation*}
and
\begin{equation*}
    \tgtwt_j \optsol_{jj} = \tgtwt_j - \tgtwt_j [1 - \srcwt_j / \tgtwt_j]^+ = \tgtwt_j - [\tgtwt_j - \srcwt_j]^+ = \srcwt_j - [\srcwt_j - \tgtwt_j]^+.
\end{equation*}

Now first,
\begin{align*}
    \sum_{j=1}^n \optsol_{ij}
    &= \optsol_{ii} + \sum_{j \ne i} \optsol_{ij} \\
    &= 1 - [1 - \srcwt_i / \tgtwt_i]^+ + \sum_{j \ne i} \frac{[1 - \srcwt_i / \tgtwt_i]^+ [\srcwt_j - \tgtwt_j]^+}
                                                             {\sum_{k=1}^n [\srcwt_k - \tgtwt_k]^+} \\
    &= 1 - \frac{[\tgtwt_i - \srcwt_i]^+ [\srcwt_j - \tgtwt_j]^+}
                {\tgtwt_i \sum_{k=1}^n [\srcwt_k - \tgtwt_k]^+} \\
    &= 1.
\end{align*}
And second,
\begin{align*}
    \sum_{i=1}^n \tgtwt_i \optsol_{ij}
    &= \tgtwt_j \optsol_{jj} + \sum_{i \ne j} \tgtwt_i \optsol_{ij} \\
    &= \srcwt_j - [\srcwt_j - \tgtwt_j]^+ + \sum_{i \ne j} \frac{\tgtwt_i [1 - \srcwt_i / \tgtwt_i]^+ [\srcwt_j - \tgtwt_j]^+}{\sum_{k=1}^n [\srcwt_k - \tgtwt_k]^+} \\
    &= \srcwt_j - \frac{[\srcwt_j - \tgtwt_j]^+ [\tgtwt_j - \srcwt_j]^+}
                    {\sum_{k=1}^n [\srcwt_k - \tgtwt_k]^+} \\
    &= \srcwt_j.
\end{align*}
Thus, $\Optsol \bs{1} = \bs{1}$ and $\Optsol^\top \tgtwtvec  = \srcwtvec$. Finally, $\Optsol \ge 0$ is evident from the definition of $\optsol_{ij}$. Hence, $\Optsol \in \Solset$.
\end{proof}

\begin{lem} \label{lem:dual_feasible}
For any two $n$-dimensional stochastic vectors ${\srcwtvec > \bs{0}}$ and ${\tgtwtvec > \bs{0}}$, the corresponding vectors $\optdualavec$ and $\optdualbvec$ and matrix $\optdualcmat$ are feasible in \eqref{eq:D}.
\end{lem}
\begin{proof}
One may verify that $\optdualavec, \optdualbvec$, and $\optdualcmat$ comprise a feasible solution to \eqref{eq:D} from Table~\ref{tab:dual_solution}.
\end{proof}

\begin{lem} \label{lem:complementary_slackness}
For any two $n$-dimensional stochastic vectors ${\srcwtvec > \bs{0}}$ and ${\tgtwtvec > \bs{0}}$, the corresponding matrices $\Optsol$ and $\optdualcmat$ satisfy the complementary slackness condition $\tr(\Sol \dualcmat^\top) = 0$.
\end{lem}
\begin{proof}
Observe that $\optdualc_{ij} = 0$ only when $i = j$ or when $i \ne j, \srcwt_i < \tgtwt_i$, and $\srcwt_j \ge \tgtwt_j$. In any other case, $\optsol_{ij} = 0$. As such $\optdualc_{ij} \optsol_{ij} = 0$ for all $i,j=1,\ldots,n$.
\end{proof}

\begin{thm} \label{thm:optimal}
For any two $n$-dimensional stochastic vectors ${\srcwtvec > \bs{0}}$ and ${\tgtwtvec > \bs{0}}$, the corresponding $\Optsol$ is optimal in \eqref{eq:P} and the corresponding $\optdualavec$, $\optdualbvec$, and $\optdualcmat$ are optimal in \eqref{eq:D}.
\end{thm}
\begin{proof}
The KKT conditions in Table~\ref{tab:KKT} are necessary and sufficient for LP optimality. For $\Optsol$, $\optdualavec$, $\optdualbvec$, and $\optdualcmat$, Lemma~\ref{lem:primal_feasible} proves primal feasibility, Lemma~\ref{lem:dual_feasible} proves dual feasibility and stationarity, and Lemma~\ref{lem:complementary_slackness} proves complementary slackness.
\end{proof}

It follows, of course, that $\tr(\Optsol) = \optdualavec^\top \srcwtvec + \optdualbvec^\top \bs{1}$. No matter, showing how this equality holds shines light on the relationship between $[\,\cdot\,]^+$ in the construction of the primal optimal solution and $\mathbb{I}_{[\,\cdot\,]}$ in the construction of the dual optimal solution:
\begin{align*}
    \optdualavec^\top \srcwtvec + \optdualbvec^\top \bs{1}
    &= \sum_{i=1}^n \left( \srcwt_i / \tgtwt_i \cdot \mathbb{I}_{[\srcwt_i < \tgtwt_i]} + \mathbb{I}_{[\srcwt_i \ge \tgtwt_i]} \right) \\
    &= \sum_{i=1}^n \left( \srcwt_i / \tgtwt_i \cdot (1 - \mathbb{I}_{[\srcwt_i \ge \tgtwt_i]}) + \mathbb{I}_{[\srcwt_i \ge \tgtwt_i]} \right) \\
    &= \sum_{i=1}^n \left( (1 - \srcwt_i / \tgtwt_i) \mathbb{I}_{[\srcwt_i \ge \tgtwt_i]} + \srcwt_i / \tgtwt_i \right) \\
    &= \sum_{i=1}^n \left( 1 - [1 - \srcwt_i / \tgtwt_i]^+ \right) \\
    &= \sum_{i=1}^n \optsol_{ii} \\
    &= \tr(\Optsol).
\end{align*}

To conclude this section, we discuss conditions under which $\Optsol$ is uniquely optimal in \eqref{eq:P}.

\begin{lem} \label{lem:opt_scaled_basis_vectors}
Every optimal solution to \eqref{eq:P} possesses
\begin{equation*}
    \sol_{ii} = \min\{p_i / q_i, 1\} = 1 - [1 - p_i / q_i]^+, \quad \forall i=1,\ldots,n.
\end{equation*}
Moreover, if an optimal solution possesses $\sol_{ii} = 1$ then its $i\textsuperscript{th}$ row is $\bs{e}_i^\top$, and if it possesses $\sol_{ii} = \frac{p_i}{q_i}$ then its $i\textsuperscript{th}$ column is $\frac{p_i}{q_i} \bs{e}_i$.
\end{lem}
\begin{proof}
For each $i=1,\ldots,n$,
\begin{align*}
    \Sol \bs{1} = \bs{1}, \Sol \ge 0 &\implies \sol_{ii} \le 1, \\
    \Sol^\top \bs{q} = \bs{p}, \Sol \ge 0 &\implies \sol_{ii} \le p_i / q_i.
\end{align*}
As such, the constraints of \eqref{eq:P} imply for each $i=1,\ldots,n$ that ${u_{ii} \le \min\{p_i / q_i, 1\} = 1 - [1 - p_i / q_i]^+}$ independently. The objective of \eqref{eq:P} is to maximize $\sum_{i=1}^n \sol_{ii}$, and $\Optsol$ does so by satisfying all of these constraints at equality. Hence, every optimal solution to \eqref{eq:P} must satisfy these constraints at equality.

If $u_{ii} = 1$, then $\Sol \bs{1} = \bs{1}, \Sol \ge 0$ imply that $u_{ij} = 0$ for all $j \ne i$. Similarly, if $u_{ii} = p_i / q_i$, then $\Sol^\top \bs{q} = \bs{p}, \Sol \ge 0$ imply that $u_{ji} = 0$ for all $j \ne i$.
\end{proof}

The following lemma, which provides a necessary and sufficient condition for the uniqueness of an optimal LP solution, comes from \citet{Mangasarian1979}.
\begin{lem} \label{lem:mangasarian}
Let $\overline{\bs{x}}$ be an optimal solution to the linear program
\begin{equation*}
    \min\{\bs{p}^\top \bs{x} : A \bs{x} = \bs{b}, C \bs{x} \ge d\}.
\end{equation*}
Let $\bs{c}_i^\top$ be the $i\textsuperscript{th}$ row of $C$, and let $C_\mathcal{E}$ be the matrix whose rows are all the $\bs{c}_i$ satisfying $\bs{c}_i^\top \overline{\bs{x}} = d_i$. Then $\overline{\bs{x}}$ is uniquely optimal if and only if
\begin{equation*}
    \widetilde{\mathcal{X}} = \{
        \tilde{\bs{x}} :
        \bs{p}^\top \tilde{\bs{x}} \le 0,
        \tilde{\bs{x}} \ne \bs{0},
        A \tilde{\bs{x}} = \bs{0},
        C_\mathcal{E} \tilde{\bs{x}} \ge \bs{0}
    \} = \varnothing.
\end{equation*}
\end{lem}

Lemma~\ref{lem:mangasarian} essentially states that for any $\tilde{\bs{x}} \in \widetilde{X}$ and sufficiently small $\epsilon > 0$, $\overline{\bs{x}}$ and $\overline{\bs{x}} + \epsilon \tilde{\bs{x}}$ are both feasible and optimal solutions. With this in mind, we consider the conditions under which
\begin{align*}
    \Orthset \triangleq \{
        \Optorth \in \mathbb{R}^{n \times n} :
        \tr(\Optorth) \ge 0,
        \Optorth \ne 0,
        \Optorth^\top \tgtwtvec = \bs{0},
        \Optorth \bs{1} = \bs{0}, & \\
        \optorth{ij} \ge 0, \forall i,j=1,\ldots,n~\text{s.t.}~\optsol_{ij} = 0
    &\} = \varnothing.
\end{align*}

\begin{thm} \label{thm:unique_cond}
For any two $n$-dimensional stochastic vectors ${\srcwtvec > \bs{0}}$ and ${\tgtwtvec > \bs{0}}$, define
\begin{align*}
    \mathcal{G} &= \{i : p_i > q_i\}, \\
    \mathcal{E} &= \{i : p_i = q_i\}, \\
    \mathcal{L} &= \{i : p_i < q_i\}.
\end{align*}
The corresponding matrix $\Optsol \in \mathbb{R}^{n \times n}$ is uniquely optimal in \eqref{eq:P} if and only if $|\mathcal{G}| \le 1$ or $|\mathcal{L}| \le 1$.
\end{thm}
\begin{proof}
We first consider the case of $|\mathcal{G}| \le 1$. Note that $|\mathcal{G}| = 0$ if and only if $\srcwtvec = \tgtwtvec$, and it follows from Lemma~\ref{lem:opt_scaled_basis_vectors} that $\Optsol = I$ is uniquely optimal. If $|\mathcal{G}| = 1$, then $|\mathcal{E} \cup \mathcal{L}| = n - 1$. In this case, $n-1$ columns of an optimal to solution to \eqref{eq:P} are determined by Lemma~\ref{lem:opt_scaled_basis_vectors}, and the remaining column is determined by $\Sol \bs{1} = \bs{1}$. Hence, $\Optsol$ is uniquely optimal. For $|\mathcal{L}| \le 1$, a similar argument applies to the rows of an optimal solution.

Now if neither $|\mathcal{G}| \le 1$ nor $|\mathcal{L}| \le 1$, then ${|\mathcal{G}| \ge 2}$ and ${|\mathcal{L}| \ge 2}$. Let $i, i' \in \mathcal{L}, i \ne i'$, let $j, j' \in \mathcal{G}, j \ne j'$, and let $\Optorth$ be a matrix with ${\optorth_{ij} = q_{i'}}$, ${\optorth_{ij'} = -q_{i'}}$, ${\optorth_{i'j} = -q_{i}}$, ${\optorth_{i'j'} = q_{i}}$, and all other entries zero.
Note that $\mathcal{G} \cap \mathcal{L} = \varnothing$, so $\tr(\Optorth) = 0$.
The matrix $\Optorth$ trivially satisfies $\Optorth \ne 0$, $\Optorth^\top \tgtwtvec = \bs{0}$, and $\Optorth \bs{1} = \bs{0}$.
Finally, owing to the fact that $p_i < q_i$, $p_{i'} < q_{i'}$, $p_j > q_j$, and $p_{j'} > q_{j'}$, all of $\optsol_{ij}$, $\optsol_{ij'}$, $\optsol_{i'j}$, and $\optsol_{i'j'}$ are strictly positive. Accordingly, the corresponding entries of $\Optorth$ are not subject to the only other constraint of $\Orthset$.
Hence, $\Optorth \in \Orthset \ne \varnothing$.
Following Lemma~\ref{lem:mangasarian}, $\Optsol$ is not uniquely optimal in \eqref{eq:P}.
\end{proof}

\begin{cor}
For any two $n$-dimensional stochastic vectors ${\srcwtvec > \bs{0}}$ and ${\tgtwtvec > \bs{0}}$, the corresponding $\Optsol$ is uniquely optimal in \eqref{eq:P} if $n=1, 2, 3$.
\end{cor}

\begin{example}
Let $\srcwtvec^\top = [\sfrac{1}{2}, \sfrac{1}{4}, \sfrac{1}{8}, \sfrac{1}{8}]$ and $\tgtwtvec^\top = [\sfrac{1}{20}, \sfrac{1}{5}, \sfrac{7}{20}, \sfrac{2}{5}]$. Then $\mathcal{G} = \{1, 2\}$, $\mathcal{E} = \varnothing$, $\mathcal{L} = \{3, 4\}$, and
\begin{equation*}
    \Optsol = \left[ \begin{array}{cccc}
        1 & 0 & 0 & 0 \\
        0 & 1 & 0 & 0 \\
        \sfrac{3}{8} & \sfrac{1}{8} & \sfrac{1}{2} & 0 \\
        \sfrac{9}{16} & \sfrac{3}{16} & 0 & \sfrac{1}{4}
    \end{array} \right].
\end{equation*}
Using the construction in the proof for Theorem~\ref{thm:unique_cond}, we obtain
\begin{equation*}
    \Optorth = \left[ \begin{array}{cccc}
        0 & 0 & 0 & 0 \\
        0 & 0 & 0 & 0 \\
        \sfrac{-1}{2} & \sfrac{1}{2} & 0 & 0 \\
        \sfrac{1}{4} & \sfrac{-1}{4} & 0 & 0
    \end{array} \right] \in \Orthset.
\end{equation*}
From this, we see that
\begin{equation*}
    \Optsol + \min\{\sfrac{\optsol_{31}}{q_4}, \sfrac{\optsol_{42}}{q_3}\} \Optorth = \left[ \begin{array}{cccc}
        1 & 0 & 0 & 0 \\
        0 & 1 & 0 & 0 \\
        0 & \sfrac{1}{2} & \sfrac{1}{2} & 0 \\
        \sfrac{3}{4} & 0 & 0 & \sfrac{1}{4}
    \end{array} \right]
\end{equation*}
is likewise feasible and optimal in \eqref{eq:P}.

\end{example}

%% file: sections/04_properties.tex
\section{Properties of $\Optsol$ and its Inverse} \label{section:properties}

In this section, we identify the matrix $\Optinv$ that is the inverse of $\Optsol$, then we discuss the properties of these matrices. First, we show that $\Optinv$ is a $Z$-matrix and moreover an $M$-matrix, and that $\Optsol$ is accordingly an inverse $M$-matrix (\textit{i.e.}, a matrix that is the inverse of an $M$-matrix). These classes of matrices are noteworthy for their numerous and varied characterizations \citep{Fiedler1962,Berman1994}. Independently, we show that the diagonal entries of $\Optsol$ and $\Optinv$ are the eigenvalues of those matrices.

For stochastic vectors $\srcwtvec > \bs{0}$ and $\tgtwtvec > \bs{0}$, let $\Optinv$ be the corresponding matrix with entries
\begin{equation}
    \optinv_{ij} \triangleq \begin{cases}
        1 - [1 - \tgtwt_i / \srcwt_i]^-,
        & i = j, \\[.5em]
        0,
        & i \ne j~\text{and}~\srcwtvec = \tgtwtvec, \\[0.25em]
        \displaystyle \frac{[1 - \tgtwt_i / \srcwt_i]^- [\srcwt_j - \tgtwt_j]^+}
                           {\sum_{k=1}^{n} [\srcwt_k - \tgtwt_k]^+},
        & i \ne j~\text{and}~\srcwtvec \ne \tgtwtvec.
    \end{cases}
    \label{eq:Vopt}
\end{equation}
We now prove that $\Optinv$ is the inverse of $\Optsol$.

\begin{thm} \label{thm:UV_inverses}
For any two $n$-dimensional stochastic vectors ${\srcwtvec > \bs{0}}$ and ${\tgtwtvec > \bs{0}}$, the corresponding matrices $\Optsol$ and $\Optinv$ satisfy $\Optsol \Optinv = I$.
\end{thm}
\begin{proof}
If $\srcwtvec = \tgtwtvec$, then $\Optsol = \Optinv = I$, and the result is trivially true. Otherwise, observe that
\begin{gather*}
    \optsol_{kk}
    = 1 - [1 - \srcwt_k / \tgtwt_k]^+
    = \min\{\srcwt_k / \tgtwt_k, 1\}, \\
    \optinv_{kk}
    = 1 - [1 - \tgtwt_k / \srcwt_k]^-
    = \max\{\tgtwt_k / \srcwt_k, 1\}.
\end{gather*}
If $\srcwt_k \ge \tgtwt_k$, then $\optsol_{kk} = 1$ and $\optinv_{kk} = 1$. If $\srcwt_k < \tgtwt_k$, then $\optsol_{kk} = \srcwt_k / \tgtwt_k$ and $\optinv_{kk} = \tgtwt_k / \srcwt_k$. Hence, $\optsol_{kk} \optinv_{kk} = 1$. Additionally, observe that $[\srcwt_k - \tgtwt_k]^+ [1 - \tgtwt_k / \srcwt_k]^- = 0$, so
\begin{equation*}
    \optsol_{ik} \optinv_{kj}
    =
    \frac{[1 - \srcwt_i / \tgtwt_i]^+ [\srcwt_k - \tgtwt_k]^+}
         {\sum_{l=1}^{n} [\srcwt_l - \tgtwt_l]^+}
    \cdot 
    \frac{[1 - \tgtwt_k / \srcwt_k]^- [\srcwt_j - \tgtwt_j]^+}
         {\sum_{l=1}^{n} [\srcwt_l - \tgtwt_l]^+}
    = 0
\end{equation*}
for all $i, j, k = 1, \ldots, n, i \ne k, k \ne j$. As such,
\begin{equation*}
    \left[ \Optsol \Optinv \right]_{ij}
    = \sum_{k=1}^n \optsol_{ik} \optinv_{kj}
    = \begin{cases}
        1~\text{if}~ i = j, \\
        0~\text{otherwise.}
    \end{cases}
\end{equation*}
\end{proof}

Since $\bs{1}$ is a fixed point of $\Optsol$, it is also a fixed point of its inverse $\Optinv$. Hence, $\Optinv \bs{1} = \bs{1}$. Also,
\begin{equation*}
    \Optsol^\top \srcwtvec = \tgtwtvec \implies \Optinv^\top \Optsol^\top \srcwtvec = \Optinv^\top \tgtwtvec \implies \Optinv^\top \tgtwtvec = \srcwtvec.
\end{equation*}
Geometrically, $\Optsol$ is the matrix representation of a simplex that has a $\srcwtvec$-weighted barycenter of $\tgtwtvec$ and is enclosed by the standard simplex. Its inverse $\Optinv$ is the matrix representation of a simplex that has a $\tgtwtvec$-weighted barycenter of $\srcwtvec$ and is an enclosure of the standard simplex.


\subsection{Classification of $\Optsol$ and $\Optinv$}

Using the fact that $\Optsol \Optinv = I$, we now classify these matrices.

\begin{lem} \label{lem:V_is_Zmatrix}
For any two $n$-dimensional stochastic vectors ${\srcwtvec > \bs{0}}$ and ${\tgtwtvec > \bs{0}}$, the corresponding $\Optinv$ is a $Z$-matrix.
\end{lem}
\begin{proof} \label{prf:Zmatrix}
A $Z$-matrix is a real-valued square matrix whose off-diagonal entries are nonpositive \citep{Berman1994}. From \eqref{eq:Vopt}, we see that the offdiagonal entries of $\Optinv$ are either zero (if $\srcwtvec = \tgtwtvec$) or nonpositive (if $\srcwtvec \ne \tgtwtvec$).
\end{proof}

\begin{thm} \label{thm:V_is_Mmatrix}
For any two $n$-dimensional stochastic vectors ${\srcwtvec > \bs{0}}$ and ${\tgtwtvec > \bs{0}}$, the corresponding $\Optinv$ is an $M$-matrix.
\end{thm}
\begin{proof} \label{prf:V_is_Mmatrix}
For a matrix $A$ to be an $M$-matrix, it is necessary and sufficient for $A$ to be a $Z$-matrix that is nonsingular with ${A^{-1} \ge 0}$ \citep{Berman1994}, which for $\Optinv$ follows from Theorem~\ref{lem:primal_feasible}, Theorem~\ref{thm:UV_inverses}, and Lemma~\ref{lem:V_is_Zmatrix}.
\end{proof}

\begin{cor}
For any two $n$-dimensional stochastic vectors ${\srcwtvec > \bs{0}}$ and ${\tgtwtvec > \bs{0}}$, the corresponding $\Optsol$ is an inverse $M$-matrix.
\end{cor}

\subsection{Eigenvalues of $\Optsol$ and $\Optinv$}

In Section~\ref{section:solutions} we showed that $\Optsol$ is an optimal solution to the matrix optimization problem \eqref{eq:P}, which we noted is an eigenvalue optimization problem. Interestingly, the eigenvalues of $\Optsol$ and $\Optinv$ are their diagonal entries, which we now prove.

\begin{lem} \label{lem:UV_scaled_basis_vectors}
For any two $n$-dimensional stochastic vectors ${\srcwtvec > \bs{0}}$ and ${\tgtwtvec > \bs{0}}$, for all $i = 1, 2, \ldots, n,$ either the $i\textsuperscript{th}$ row or $i\textsuperscript{th}$ column of the corresponding $\Optsol$ consists entirely of zeros except for $\optsol_{ii} > 0$, and either the $i\textsuperscript{th}$ row or $i\textsuperscript{th}$ column of the corresponding $\Optinv$ consists entirely of zeros except for $\optinv_{ii} > 0$.
\end{lem}
\begin{proof}
If $\srcwtvec = \tgtwtvec$, then $\Optsol = \Optinv = I$. In such a case, both the $i\textsuperscript{th}$ row and $i\textsuperscript{th}$ column of $\Optsol$ are all zeros except for $\optsol_{ii} = 1$. Similarly, both the $i\textsuperscript{th}$ row and $i\textsuperscript{th}$ column of $\Optinv$ are all zeros except for $\optinv_{ii} = 1$.

Otherwise, if $\srcwtvec \ne \tgtwtvec$, recall for $i \ne j$ that
\begin{align*}
    \optsol_{ij} =
    \frac{[1 - \srcwt_i / \tgtwt_i]^+ [\srcwt_j - \tgtwt_j]^+}
         {\sum_{k=1}^{n} [\srcwt_k - \tgtwt_k]^+}, \\
    \optinv_{ij} =
    \frac{[1 - \tgtwt_i / \srcwt_i]^- [\srcwt_j - \tgtwt_j]^+}
         {\sum_{k=1}^{n} [\srcwt_k - \tgtwt_k]^+}.
\end{align*}
Now observe that for all $i,j=1,\ldots,n, i \ne j$, either
\begin{equation*}
    \srcwt_i \ge \tgtwt_i \implies [1 - \srcwt_i / \tgtwt_i]^+ = [1 - \tgtwt_i / \srcwt_i]^- = 0 \implies \optsol_{ij} = \optinv_{ij} = 0,
\end{equation*}
or
\begin{equation*}
    \srcwt_i \le \tgtwt_i \implies [\srcwt_i - \tgtwt_i]^+ \implies \optsol_{ji} = \optinv_{ji} = 0.
\end{equation*}
\end{proof}

\begin{thm} \label{thm:UV_eigenvalues}
For any two $n$-dimensional stochastic vectors ${\srcwtvec > \bs{0}}$ and ${\tgtwtvec > \bs{0}}$, the eigenvalues of the corresponding $\Optsol$ are $\optsol_{11}, \optsol_{22}, \ldots, \optsol_{nn}$, and the eigenvalues of the corresponding $\Optinv$ are similarly $\optinv_{11}, \optinv_{22}, \ldots, \optinv_{nn}$.
\end{thm}
\begin{proof}
From Lemma~\ref{lem:UV_scaled_basis_vectors}, for $i=1,\ldots,n$ either the $i\textsuperscript{th}$ row or $i\textsuperscript{th}$ column of $\Optsol$ consists of all zeros except for $\optsol_{ii} > 0$. The same is true of $\Optinv$. Using Laplace expansion, one may thus expand strategically on rows or columns to get characteristic equations
\begin{gather*}
    \det(\Optsol - \lambda I) = \prod_{i=1}^n (\optsol_{ii} - \lambda), \\
    \det(\Optinv - \lambda I) = \prod_{i=1}^n (\optinv_{ii} - \lambda).
\end{gather*}
Ergo, the eigenvalues of $\Optsol$ are $\optsol_{11}, \optsol_{22}, \ldots, \optsol_{nn}$ and, similarly, the eigenvalues of $\Optinv$ are $\optinv_{11}, \optinv_{22}, \ldots, \optinv_{nn}$.
\end{proof}

While $\Optsol$, among $\Sol \in \Solset$, maximizes $\tr(\Sol)$, it does not necessarily maximize $|\det(\Sol)|$. That is, among $\Sol \in \Solset$, the matrix $\Optsol$ does not necessarily correspond to a simplex of maximum unsigned volume. The reason is that trace maximization inherently promotes a simplex orientation that sites the $i\textsuperscript{th}$ vertex near to $\bs{e}_i$. Consider the following relaxation of \eqref{eq:P} involving the joint optimization of $\Sol \in \Solset$ and $\Perm \in \Permset$, the set of all $n \times n$ permutation matrices:
\begin{equation}
    \max\{\tr(\Perm \Sol) : \Sol \in \Solset, \Perm \in \Permset\}. \label{eq:P-Perm} \tag{P-Perm}
\end{equation}
The goal of this problem is to identify a $\Sol \in \Solset$ whose premultiplication with an appropriate $\Perm \in \Permset$ yields a trace-maximizing matrix. This relaxation of \eqref{eq:P} promotes no particular simplex orientation. Rather, it allows the simplex associated with $\Sol$ to be oriented freely via $\Perm \in \Permset$.

For fixed $\Perm \in \Permset$, a matrix $\Sol \in \Solset$ that maximizes $\tr(\Perm \Sol)$ is the $\Optsol$ corresponding to stochastic vectors $\srcwtvec > 0$ and $M^\top \tgtwtvec > 0$. An optimal $\Perm$ to \eqref{eq:P-Perm} may thus be identified via
\begin{maxi*}
    {}{\sum_{i \in \mathcal{R}} \sum_{j \in \mathcal{R}} (1 - [1 - \srcwt_j / \tgtwt_i]^+) \perm_{ij}}{}{}
    \addConstraint{\Perm \bs{1} = \bs{1}}
    \addConstraint{\Perm^\top \bs{1} = \bs{1}}
    \addConstraint{\Perm \in \{0, 1\}^{n \times n}}
\end{maxi*}
This assignment problem features discrete variables. However, the variable coefficients in the constraints comprise a totally unimodular matrix and the constant terms in the constraints are all integral. As such, every optimal extreme point solution to the linear programming relaxation is an optimal $\Perm$ to \eqref{eq:P-Perm} \citep{Wolsey2020}.



That \eqref{eq:P-Perm} is a relaxation of \eqref{eq:P} means there may exist a $\widehat{\Sol} \in \Solset$ and $\widehat{\Perm} \in \Permset$ for which $\tr(\widehat{\Perm}\widehat{\Sol}) \ge \tr(\Optsol)$. Of course, if $\tr(\widehat{\Perm}\widehat{\Sol}) > \tr(\Optsol)$, it is possible that $\det(\widehat{\Perm}\widehat{\Sol}) > \det(\Optsol)$ and, since the absolute value of the determinant of a matrix is invariant under permutation, $|\det(\widehat{\Sol})| > |\det(\Optsol)|$. We demonstrate this now in an example.

\begin{example} \label{ex:3}
Let $\srcwtvec^\top = [\sfrac{3}{10}, \sfrac{2}{5}, \sfrac{1}{10}, \sfrac{1}{5}]$ and $\tgtwtvec^\top = [\sfrac{1}{8}, \sfrac{3}{8}, \sfrac{3}{10}, \sfrac{1}{5}]$. Then
\begin{equation*}
    \Optsol = \left[ \begin{array}{cccc}
        1 & 0 & 0 & 0 \\
        0 & 1 & 0 & 0 \\
        \sfrac{7}{12} & \sfrac{1}{12} & \sfrac{1}{3} & 0 \\
        0 & 0 & 0 & 1
    \end{array} \right]
\end{equation*}
is the optimal solution to \eqref{eq:P}, and
\begin{equation*}
    \widehat{\Perm} = \left[ \begin{array}{cccc}
        0 & 0 & 1 & 0 \\
        0 & 1 & 0 & 0 \\
        1 & 0 & 0 & 0 \\
        0 & 0 & 0 & 1
    \end{array} \right],~
    \widehat{\Sol} = \left[ \begin{array}{cccc}
        0 & \sfrac{1}{5} & \sfrac{4}{5} & 0 \\
        0 & 1 & 0 & 0 \\
        1 & 0 & 0 & 0 \\
        0 & 0 & 0 & 1
    \end{array} \right]
\end{equation*}
are the optimal solution to \eqref{eq:P-Perm}.
For this example, we have $\tr(\widehat{\Perm} \widehat{\Sol}) = \sfrac{19}{5} > \sfrac{10}{3} = \tr(\Optsol)$
and $|\det(\widehat{\Sol})| = \sfrac{4}{5} > \sfrac{1}{3} = |\det(\Optsol)|$.
The simplices associated with $\Optsol$ and $\widehat{\Sol}$ are illustrated in Figure~\ref{fig:example3}. Indeed, the simplex associated with $\widehat{\Sol}$ has a visibly larger volume than that associated with $\Optsol$.
\begin{figure}
    \centering
    \begin{minipage}{.45\linewidth}
        \includegraphics[trim={55pt 60pt 75pt 80pt},clip,width=\linewidth]{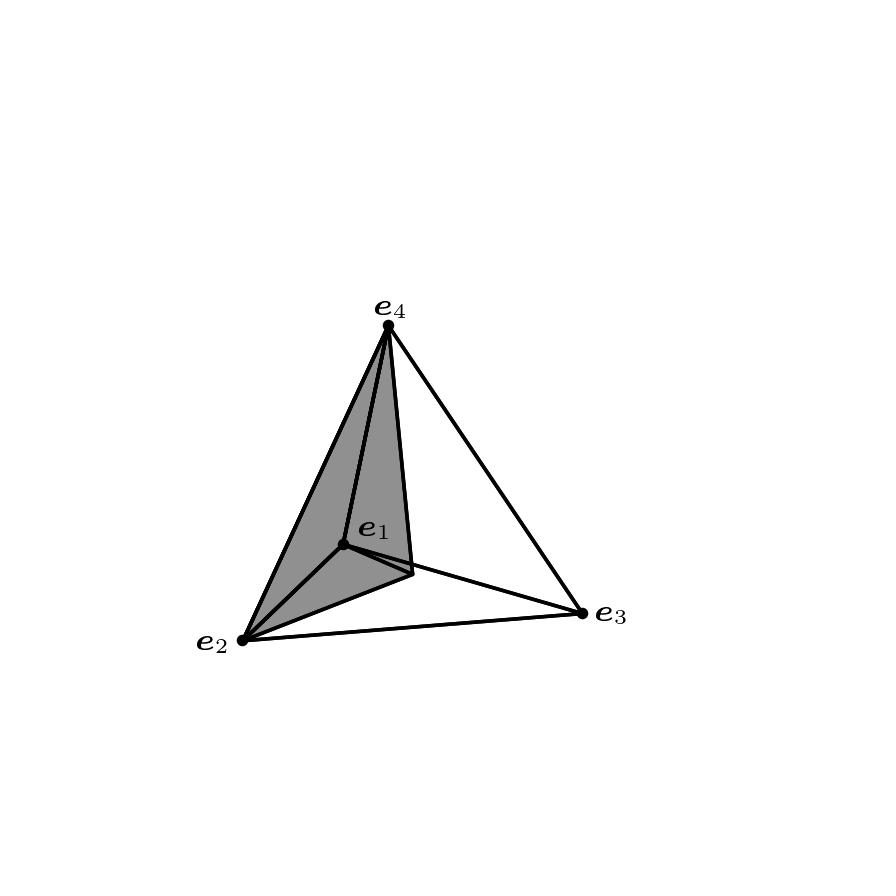}
    \end{minipage}
    \hfill
    \begin{minipage}{.45\linewidth}
        \includegraphics[trim={55pt 60pt 75pt 80pt},clip,width=\linewidth]{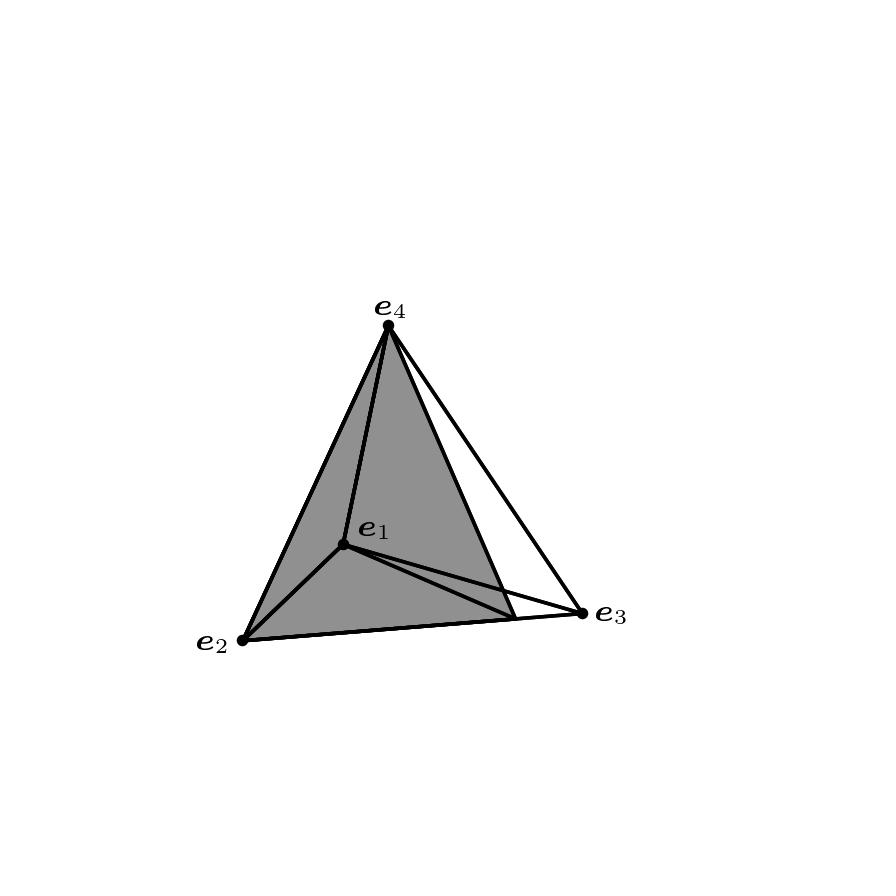}
    \end{minipage}
    \caption{The simplices associated with $\Optsol$ (left) and $\widehat{\Sol}$ (right) from Example~\ref{ex:3}.}
    \label{fig:example3}
\end{figure}
\end{example}


%% file: sections/05_application.tex
\section{Application} \label{section:application}

Let $X$ be a random variable (RV) with support $\Omega$ and associated CDF $F(x)$. We wish to synthesize conditional CDFs with respect to a discrete RV $Y$ with support $\mathcal{N} = \{1, 2, \ldots, n\}$ according to
\begin{equation}
	F(x\,|\,Y = i) = \srctm_i(F(x)), \quad \forall i \in \mathcal{N}. \label{eq:synth_cond_CDF}
\end{equation}
Now define $\srcwt_i = \prob(Y = i)$. Also, define $\srcwtvec = [\srcwt_1, \srcwt_2, \ldots, \srcwt_n]^\top$, and $\srctmvec(z) = [\srctm_1(z), \srctm_2(z), \ldots, \srctm_n(z)]^\top$. For such synthesized conditional CDFs to be coherent, $\srcwtvec$ and $\srctmvec(z)$ must satisfy
\begin{align}
    & \srcwtvec^\top \ge \bs{0}, \label{eq:wt_nonnegative} \\
    & \srcwtvec^\top \bs{1} = 1, \label{eq:wt_sum_to_one} && \\
	& \srctmvec(0) = \bs{0}, \label{eq:tm_0} \\
    & \srctmvec(1) = \bs{1}, \label{eq:tm_1} \\
    & \srctmvec(z) \le \srctmvec(\hat{z}), &&  \forall z, \hat{z} \in [0, 1], z \le \hat{z}, \label{eq:tm_monotone} \\
    & \srcwtvec^\top \srctmvec(z) = z, && \forall z \in [0, 1]. \label{eq:wt_tm_lotop}
\end{align}
Conditions \eqref{eq:wt_nonnegative} and \eqref{eq:wt_sum_to_one} are simply the first two axioms of probability. Conditions \eqref{eq:tm_0}-\eqref{eq:tm_monotone} ensure for each $i \in \mathcal{N}$ that $\srctm_i(F(x))$ is the CDF of an RV with support $\Omega_i \subseteq \Omega$. Condition \eqref{eq:wt_tm_lotop}, equivalent to $\srcwtvec^\top \srctmvec(F(x)) = F(x), \forall x \in \Omega$, is the law of total probability.

Suppose there exist independent RVs $X_1, X_2, \ldots, X_m$ all identically distributed according to $F(x)$. 
Essentially, we wish to partition the RVs into sets enumerated $i=1,\ldots,n$ with relative sizes $\srcwt_1, \ldots, \srcwt_n$, and then synthesize conditional CDFs to force a distinction between RVs from different sets while retaining $F(x)$ as the marginal CDF.

As we showed in our original work on this application \citep{Austgen2022}, and as we show in the sequel, one may identify with little effort special cases of $\srcwtvec$ and $\srctm(z)$ that satisfy \eqref{eq:wt_nonnegative}-\eqref{eq:wt_tm_lotop}. However, the question that motivates $\Solset$ and $\Optsol$ pertains to more general cases. Specifically, for any given $\srcwtvec$ that satisfies \eqref{eq:wt_nonnegative} and \eqref{eq:wt_sum_to_one}, do there exist unique $\srctm_i(z), \forall i \in \mathcal{N}$ that satisfy \eqref{eq:tm_0}-\eqref{eq:wt_tm_lotop}? If so, which exhibit maximum distinguishability (a notion we discuss further after the following theorem)?

\begin{thm} \label{thm:tm_existence}
Suppose $\srcwtvec > \bs{0}$ and $\srctmvec(z)$ satisfy \eqref{eq:wt_nonnegative}-\eqref{eq:wt_tm_lotop}, and let $\tgtwtvec > \bs{0}$ be a given vector that satisfies \eqref{eq:wt_nonnegative} and \eqref{eq:wt_sum_to_one}. Then for any $\Sol \in \Solset$, $\tgtwtvec$ and $\tgttmvec(z) = \Sol \srctmvec(z)$ satisfy \eqref{eq:tm_0}-\eqref{eq:wt_tm_lotop}.
\end{thm}
\begin{proof}
Given $\srctmvec(0) = \bs{0}$, it follows that ${\tgttmvec(0) = \Sol \srctmvec(0) = \bs{0}}$. Given $\srctmvec(1) = \bs{1}$ and $\Sol \bs{1} = \bs{1}$, it follows that ${\tgttmvec(1) = \Sol \srctmvec(1) = \Sol \bs{1} = \bs{1}}$. That $\Sol \bs{1} = \bs{1}$ and $\Sol \ge 0$ means that each $\tgttm_i(z)$ is a convex combination of functions that are nondecreasing on $z \in [0, 1]$. Each $\tgttm_i(z)$ is, as such, also nondecreasing on $z \in [0, 1]$. Finally, recall that ${\tgtwtvec^\top \Sol = \srcwtvec^\top}$. Given $\srcwtvec^\top \srctmvec(z) = z$, it follows that ${\tgtwtvec^\top \tgttmvec(z) = \tgtwtvec^\top \Sol \srctmvec(z) = \srcwtvec^\top \srctmvec(z) = z}$.
\end{proof}

This theorem guarantees that for any given $\tgtwtvec$ that satisfies \eqref{eq:wt_nonnegative} and \eqref{eq:wt_sum_to_one} there exists a $\tgttmvec(z)$ that satisfies \eqref{eq:tm_0}-\eqref{eq:wt_tm_lotop} provided ${\Solset \ne \varnothing}$ and there exist $\srcwtvec$ and $\srctmvec(z)$ that satisfy \eqref{eq:wt_nonnegative}-\eqref{eq:wt_tm_lotop}. Of course, $\Optsol \in \Solset \implies \Solset \ne \varnothing$. Additionally, note that $\Optsol$ maximizes the contribution of $\srctm_i(z)$ to $\tgttm_i(z)$ in the mapping $\tgttmvec(z) = U \srctmvec(z)$. So if $\srctm_i(z), i=1,\ldots,n$ are distinguishable in the sense that $\srctm_i(z) \ge \srctm_{i+1}, \forall z \in [0, 1]$ and $\int_0^1 \left( \srctm_i(z) - \srctm_{i+1}(z) \right) \diff z > 0$ for all $1 \le i < n$, mapping $\srctmvec(z)$ to $\tgttmvec(z)$ via $\Optsol$ in a sense maximizes the distinguishability of $\tgttm_i(z), i=1,\ldots,n$. Regarding the existence of $\srcwtvec$ and $\srctmvec(z)$ that satisfy \eqref{eq:wt_nonnegative}-\eqref{eq:wt_tm_lotop}, we now present a pair of examples.

First, for $i=1,\ldots,n$ let $\srcwt_i = 1 / n$ such that $\srcwtvec$ satisfies \eqref{eq:wt_nonnegative} and \eqref{eq:wt_sum_to_one}, and let $\srctm_i(z)$ be the CDF of the continuous uniform distribution with support interval $[\frac{i-1}{n}, \frac{i}{n})$:
\begin{equation} \label{eq:tm_uniform}
    \srctm_i(z) = \begin{cases}
        0, & z < \frac{i-1}{n}, \\
        n z - (i - 1), & \frac{i-1}{n} \le z < \frac{i}{n}, \\
        1, & \frac{i}{n} \le z.
    \end{cases}
\end{equation}
Since each $\srctm_i(z)$ is a CDF, $\srctmvec(z)$ satisfies \eqref{eq:tm_0}-\eqref{eq:tm_monotone}. For any $\hat{z} \in [0, 1)$, ${\exists !~\hat{i} \in \{1, 2, \ldots, n\}}$ such that $\hat{z} \in [\frac{\hat{i}-1}{n}, \frac{\hat{i}}{n})$. So ${\srctm_i(\hat{z}) = 1}$ for ${i \in \{1, \ldots, \hat{i} - 1\}}$, ${\srctm_{\hat{i}}(\hat{z}) = n \hat{z} - (\hat{i} - 1)}$, and ${\srctm_i(\hat{z}) = 0}$ for ${i \in \{\hat{i} + 1, \ldots, n\}}$. Ergo,
\begin{equation*}
    \srcwtvec^\top \srctmvec(z)
    = \frac{1}{n} \left[ (\hat{i} - 1) \cdot 1 + n z - (\hat{i} - 1) + (n - \hat{i}) \cdot 0 \right]
    = z,
\end{equation*}
so $\srcwtvec$ and $\srctmvec(z)$ also satisfy \eqref{eq:wt_tm_lotop}.

Clearly, ${\srctm_1(z) \ge \srctm_2(z) \ge \cdots \ge \srctm_n(z)}$. The 1-Wasserstein distance measured between consecutive transformations is
\begin{align*}
    \int_0^1 \!\left( \srctm_i(z) - \srctm_{i+1}(z) \right) \diff z
    &= \int_{\frac{i-1}{n}}^{\frac{i}{n}} \!\left( n z - (i - 1) \right) \diff z
       + \int_{\frac{i}{n}}^{\frac{i+1}{n}} \!\left( n z - i \right) \diff z \\
    &= n \int_{\frac{i-1}{n}}^{\frac{i}{n}} \!\left( z - \tfrac{i - 1}{n} \right) \diff z
       + n \int_{\frac{i}{n}}^{\frac{i+1}{n}} \!\left( z - \tfrac{i}{n} \right) \diff z \\
    &= 2n \int_{0}^{\frac{1}{n}} z \diff z \\
    &= \frac{1}{n}.
\end{align*}
That is, $\srctm_1(z), \srctm_2(z), \ldots, \srctm_n(z)$ are equally spaced according to the $1$-Wasserstein distance.

Second, we consider $\srctm_i(z)$ based on the CDFs of beta distributions. Note that the probability density function (PDF) for a beta distribution with integer parameters $i$ and $n - i + 1$ is
\begin{align*}
    f(z; i, n-i+1)
    = \!\binom{n}{i} i z^{i-1} (1-z)^{n-i}
    = \!\sum_{j=0}^{n-i} \binom{n}{i} \binom{n-i}{j} i (-1)^j z^{i+j-1},
\end{align*}
For $i=1,\ldots,n$ again let $\srcwt_i = 1 / n$, and let
\begin{equation} \label{eq:tm_beta}
    \srctm_i(z) = F(z; i, n-i+1)
    = \sum_{j=0}^{n-1} \binom{n}{i} \binom{n-i}{j} \left( \frac{i}{i+j} \right) (-1)^j z^{i+j}.
\end{equation}
Again, since each $\srctm_i(z)$ is a CDF, $\srctmvec(z)$ satisfies \eqref{eq:tm_0}-\eqref{eq:tm_monotone}. Also,
\begin{align*}
\srcwtvec^\top \srctmvec(z)
&= \frac{1}{n} \sum_{i=1}^n \sum_{j=0}^{n-i} \binom{n}{i} \binom{n-i}{j} \left( \frac{i}{i+j} \right) (-1)^j z^{i+j} \\
&= \frac{1}{n} \sum_{i=1}^n \sum_{k=i}^{n} \binom{n}{i} \binom{n-i}{k-i} \left( \frac{i}{k} \right) (-1)^{k-i} z^k \\
&= \frac{1}{n} \sum_{k=1}^n \sum_{i=1}^k \binom{n}{i} \binom{n-i}{k-i} \left(\frac{i}{k}\right) (-1)^{k-i} z^k \\
&= \frac{1}{n} \sum_{k=1}^n \left( \sum_{i=1}^k \binom{k-1}{i-1} (-1)^{-i} \right) \binom{n}{k} (-1)^k z^k \\
&= z
\end{align*}
Above, the final equality follows from the inner sum being an alternating sum of binomial coefficients. For $k>1$, the inner sum is zero. For $k=1$, the inner sum is $-1$. Hence, $\srcwtvec$ and $\srctmvec(z)$ satisfy \eqref{eq:wt_tm_lotop}.

Since $F(z; \alpha, \beta) \ge F(z; \alpha', \beta'), \forall z \in [0, 1]$ if and only if $\alpha \le \alpha'$ and $\beta \ge \beta'$ \citep{Lisek1978,Arab2021}, we have $\srctm_1(z) \ge \srctm_2(z) \ge \cdots \ge \srctm_n(z)$. The beta distribution is related to the Bernstein basis polynomials. Let $b_i^n(z) = \binom{n}{i} z^i (1-z)^{n-i}, i=0,\ldots,n$ be the $n+1$ degree-$n$ Bernstein basis polynomials \citep{Farouki2012}. Then, ${S_i(z) = F(z; i, n-i+1) = n \int_0^z b_{i-1}^{n-1}(z') \diff z'}$. Following properties of indefinite and definite integrals of Bernstein basis polynomials, the 1-Wasserstein distance measured between consecutive transformations is
\begin{align*}
    \int_0^1 \left( \srctm_i(z) - \srctm_{i+1}(z) \right) \diff z
    &= \int_0^1 \int_0^z \left( n b_{i-1}^{n-1}(z') - n b_{i}^{n-1}(z') \right) \diff z' \diff z \\
    &= \int_0^1 \left( \sum_{j=i}^n b_j^n(z) - \sum_{j=i+1}^{n} b_j^n(z) \right) \diff z \\
    &= \int_0^1 b_i^n(z) \diff z \\
    &= \frac{1}{n+1}.
\end{align*}
Again, $\srctm_1(z), \srctm_2(z), \ldots, \srctm_n(z)$ are equally spaced according to the $1$-Wasserstein distance.

We now conclude this section with an example that demonstrates how solutions to \eqref{eq:P} or \eqref{eq:P-Perm} are applied to synthesize conditional probability distributions as we described.

\begin{example} \label{ex:4}
Let
$\srcwtvec^\top\!=\![\sfrac{1}{4}, \sfrac{1}{4}, \sfrac{1}{4}, \sfrac{1}{4}]^\top$, and
$\tgtwtvec^\top\!=\![\sfrac{1}{3}, \sfrac{1}{9}, \sfrac{1}{2}, \sfrac{1}{18}]^\top$.
For these vectors, the matrices
\begin{equation*}
    \Optsol = \left[ \begin{array}{cccc}
        \sfrac{3}{4} & \sfrac{5}{48} & 0 & \sfrac{7}{48} \\
        0 & 1 & 0 & 0 \\
        0 & \sfrac{5}{24} & \sfrac{1}{2} & \sfrac{7}{24} \\
        0 & 0 & 0 & 1
    \end{array} \right]~\text{and}~
    \widehat{\Sol} = \left[ \begin{array}{cccc}
        0 & 0 & \sfrac{3}{4} & \sfrac{1}{4} \\
        0 & 0 & 0 & 1 \\
        \sfrac{7}{18} & \sfrac{1}{2} & 0 & \sfrac{1}{9} \\
        1 & 0 & 0 & 0
    \end{array} \right]
\end{equation*}
are optimal solutions to \eqref{eq:P} and \eqref{eq:P-Perm}, respectively. To accompany $\srcwtvec$, for $n=4$ we consider $\srctmvec^\text{UD}(z)$ based on continuous uniform distributions as defined in \eqref{eq:tm_uniform} and also $\srctmvec^\text{BD}(z)$ based on beta distributions as defined in \eqref{eq:tm_beta}. Finally, let $F(x) = \Phi(x)$, the CDF of a standard normal distribution.

We show
$\srctmvec\textsuperscript{UD}(z)$,
${\overline{\tgttmvec}\textsuperscript{UD}(z) = \Optsol \srctmvec\textsuperscript{UD}(z)}$,
${\widehat{\tgttmvec}\textsuperscript{UD}(z) = \widehat{\Sol} \srctmvec\textsuperscript{UD}(z)}$,
$\srctmvec\textsuperscript{BD}(z)$,
${\overline{\tgttmvec}\textsuperscript{BD}(z) = \Optsol \srctmvec\textsuperscript{BD}(z)}$,
${\widehat{\tgttmvec}\textsuperscript{BD}(z) = \widehat{\Sol} \srctmvec\textsuperscript{BD}(z)}$,
as well as their application to synthesizing conditional probability distributions based on $F(x)$ in Figure~\ref{fig:example4}.

The piecewise-linearity of each $\srctm_i^\text{UD}(z)$ leads to the perturbations of $F(x)$, a smooth function, being non-smooth. In contrast, that each $\srctm_i^\text{UD}(z)$ is a polynomial leads to the perturbations of $F(x)$ being smooth. Additionally, each $\srctm_i^\text{UD}(z)$ is non-increasing on portions of the unit interval whereas each $\srctm_i^\text{BD}(z)$ is strictly increasing. In the context of our application, the consequence is that the support set associated with $\srctm_i^\text{UD}(F(x))$ is $\Omega_i \subset \Omega$ whereas that associated with $\srctm_i^\text{BD}(F(x))$ is $\Omega_i = \Omega$.

Applying $\widehat{\Sol}$ leads to relatively larger 1-Wasserstein distances between the consecutive conditional CDFs, as seen on the right side of Figure~\ref{fig:example4}. However, the orientation  of $\widehat{\Sol}$ leads to those CDFs having a different order, which may not ideal for certain applications of the method. While $\Optsol$ does a more modest job of separating the conditional CDFs, its orientation maintains the ordering.

\begin{figure}
    \includegraphics[width=\linewidth]{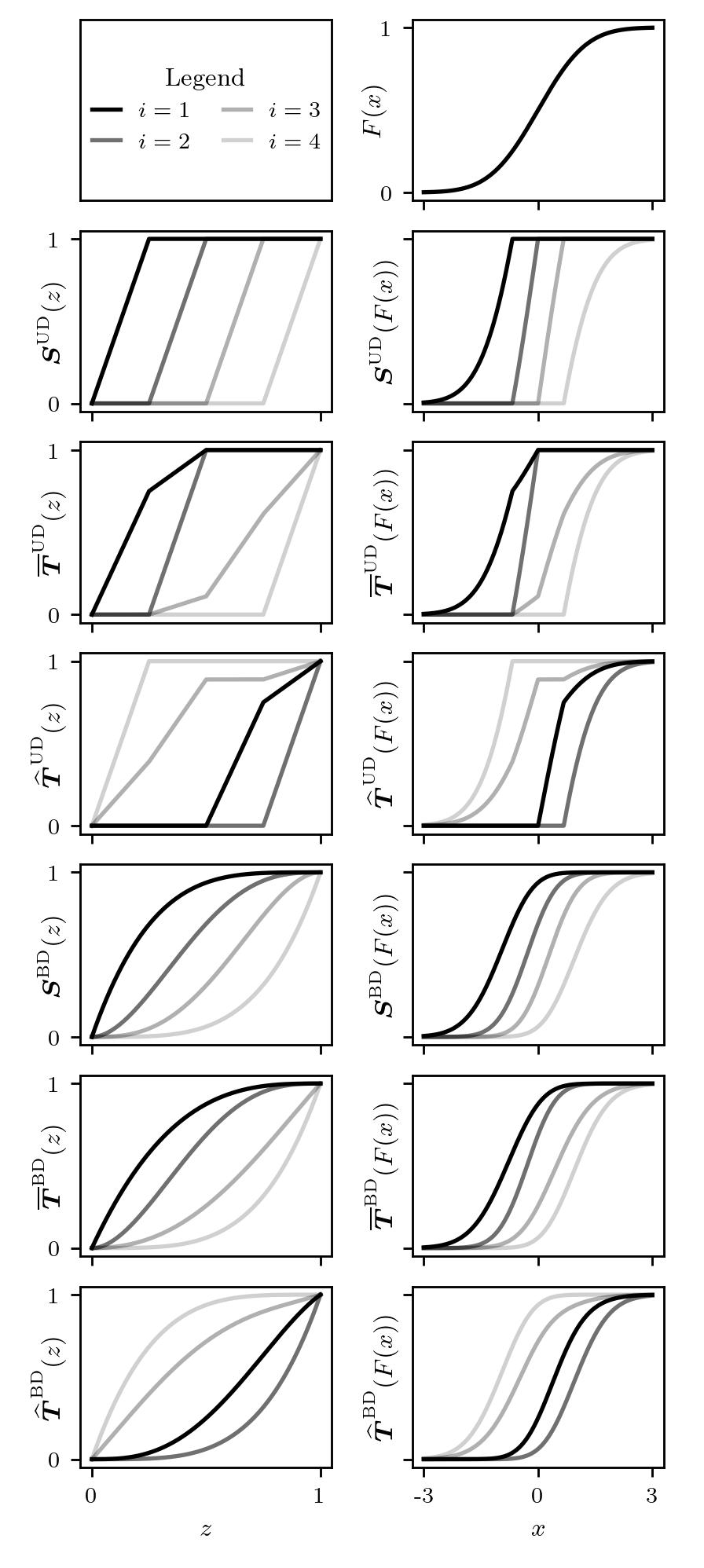}
    \caption{The transformations (left column) from Example~\ref{ex:4} and their application to synthesizing conditional CDFs based on $F(x)$ (right column).}
    \label{fig:example4}
\end{figure}

\end{example}

%% file: sections/06_conclusion.tex
\section{Conclusion} \label{section:conclusion}

We presented an optimization model involving simplices that have a given barycenter and that are enclosed by the standard simplex. We presented the analytical form of an optimal solution to the model, and the conditions under which it is the unique optimal solution. We showed the solution is an inverse $M$-matrix whose eigenvalues are the same as its diagional entries. Finally, we demonstrated how the model and its solutions apply to the task of synthesizing conditional cumulative distribution functions that, in tandem with a given discrete marginal distribution, coherently preserve a given CDF.